\documentclass[12pt]{article}
\usepackage{amsmath,amssymb,amsthm}
\usepackage[margin=3cm]{geometry}

\title{\bf The BSE property for vector-valued Frechet Lipschitz algebras}

\author{Ali Rejali$^1$\thanks{Corresponding author} \and Maryam Aghakoochaki$^2$,\thanks{2020 Mathematics Subject Classifcation. Primary: 46J05; Secondary: 46J10} }

\date{
	$^1$Department of  Pure Mathematics  Faculty of Mathematics and Statistics University of Isfahan  \\ \texttt{rejali@sci.ui.ac.ir}\\%
	$^2$University of Isfahan  \\ \texttt{mkoochaki@sci.ui.ac.ir}\\[2ex]%
	\today
}

\DeclareMathOperator{\Lip}{Lip}
\newcommand{\tnorm}[1]{{\left\vert\kern-0.25ex\left\vert\kern-0.25ex\left\vert #1 
    \right\vert\kern-0.25ex\right\vert\kern-0.25ex\right\vert}}

\theoremstyle{thmstyleone}%
\newtheorem{thm}{Theorem}

\newtheorem{pro}[thm]{Proposition}%
\newtheorem{lem}[thm]{Lemma}

\theoremstyle{thmstyletwo}%
\newtheorem{ex}{Example}%

\theoremstyle{thmstylethree}%
\newcommand{\rt}{\rightarrow}

\begin{document}

\maketitle

\begin{abstract}

  Let $(X, d)$ be a metric space with at least two elements and $(A,p_{l})$ be a commutative semisimple Frechet algebra over the scalar field $\mathbb{C}$.
The correlation between the BSE-property of the Frechet algebra $(A,p_{l})$ and $\Lip_{d}(X, A)$ is assessed. It is found and approved that if $\Lip_{d}(X, A)$ is a BSE- Frechet algebra, then so is $A$.
The opposite correlation will hold if $(A,p_{l})$ is unital.
                     
\end{abstract}
                                  \section{Introduction and preliminaries}\label{1}\setcounter{equation}{0}
The class of Frechet algebras which is an important class of locally convex algebras has been widely studied by many authors. For a full study of Frechet algebras,
 one may see (\cite{Al},\cite{G1},\cite{YH}).
A Frechet space is a metrizable complete locally convex vector space. The topology of Frechet algebra $A$ can be given by a sequence $(p_{n})$ of increasing sub- multiplicative seminorms.
 Algebra $A$ is called without order if $aA=0 $  concludes that $a=0$, $(a\in A)$.
Let $A$ be commutative and without order Frechet algebra and $\Delta(A), $ be the character space of $ A$ with the Gelfand topology. In this study, $\Delta(A)$ represents the set
of all
 non-zero multiplicative
linear functionals over $A$.
 Assume that  $C_{b}(\Delta(A))$ is the  space consisting of all complex-valued continuous and bounded functions on $\Delta(A)$.
A linear operator $T$ on $A$ is named a multiplier if $T(xy)=xT(y)$,  for  all $x,y\in A$.
The set of all multipliers on $A$ will be expressed as $M(A)$. The strong operator topology (briefly SOT- topology) on $M(A)$ is generated by the family of seminorms $\{p_{x,l}\}$ defined as
$$
p_{x,l}(T):= p_{l}(T(x))
$$
for all $x\in A$, $l\in \mathbb N$ and $T\in M(A)$.
If the Frechet algebra $A$ is  semisimple, then the Gelfand map $\Gamma: A\to\widehat A$, $f\mapsto \hat f$, is injective, or equivalently, and the following equation holds:
\begin{align*}
\underset{\varphi\in\Delta(A)}{\bigcap}ker(\varphi)=\{0\}
\end{align*}
Note that every semisimple commutative Frechet algebra is without order.  As observed in \cite{A3}, if the Frechet algebra $(A, p_{l})$ is semisimple, then 
$$(M(A),SOT)\cong (\widehat{M(A)},\mathcal{T}_{p})$$
Where $\mathcal {T}_{p}$ is pointwise topology on $\widehat{M(A)}$.

The Bochner-Schoenberg-Eberlein (BSE) is derived from the famous
theorem proved in 1980 by Bochner and Schoenberg for the group of real numbers;  \cite{SB} and \cite{P1}. The researcher in \cite{N2},  revealed that if $G$ is any locally compact abelian group, then the group algebra $L_{1}(G)$ is a BSE algebra.
 The researcher in \cite{P1},\cite{E7},\cite{E8} assessed the commutative Banach algebras that meet the Bochner-Schoenberg-Eberlein- type theorem and explained their properties.
They are introduced and assessed in    \cite{E6} the first and second types of BSE algebras. This concept is  expanded in \cite{F1} and \cite{K1}.
 
The researchers are introduced and assessed in \cite{A3}, the concept of BSE- Frechet algebra.

  The big and little Frechet $\alpha$- Lipschitz vector-valued algebra of order $\alpha$, where $\alpha\in \mathbb R$ with $\alpha>0$ was introduced in \cite{Ra}.
  The researchers are  provided a survey of the similarities and diﬀerences between Banach and
  Frechet algebras include some known results and examples. (See \cite{Al}).

 That the Lipschitz algebra $\Lip_{\alpha}(K, A)$ is a BSE-algebra if and only if $A$ is a BSE-algebra, where $K$ is a compact metric space, $A$ is a commutative unital semisimple Banach algebra, and $0<\alpha\leq 1$ is proved in  \cite{AK}.
 In this article,  this result is generalized, for any metric space $(X,d)$ and any commutative semisimple Frechet algebra $(A, p_{l})$.
That the $C_{BSE}(\Delta(\Lip_{d}(X,A)))$ can be embedded in $\Lip_{d}(X,C_{BSE}(\Delta(A))$ will be
proved in the article first, followed by proving that $\Lip_{d}(X, M(A))\underset{\rightarrow}{\subseteq} M(\Lip_{d}(X,A))$. By proving that if $\Lip_{d}(X,A)$ is a BSE- Frechet algebra, so is $A$.  If $(A, p_{l})$ is unital Frechet algbra
and  BSE- Frechet algebra, then $\Lip_{d}(X,A)$ is so,  is assessed in this article.

\setcounter{section}{0}
\section{Some basic properties of  BSE- Frechet  algebra }\label{3}\setcounter{equation}{0}
              The basic terminologies and the related information on BSE-Frechet algebras are
extracted from \cite{A3} and prove some primary, basic results, and properties related to them.

A bounded complex-valued continuous function $\sigma$ on $\Delta(A)$, is named BSE-Frechet function, if there exists a bounded set $M$ in $A$ and a positive real number $\beta_{M}$ in a sense that for every
nite complex-number $c_{1},\cdots,c_{n}$  and  the same many $\varphi_{1},\cdots,\varphi_{n}$ in $\Delta(A)$ the following inequality
$$\mid\sum_{j=1}^{n}c_{j}\sigma(\varphi_{j})\mid\leq \beta_{M}P_{M}(\sum_{j=1}^{n}c_{j}\varphi_{j})$$
holds;
where  $P_{M}$ is defined as
\begin{align*}
P_{M}(f):= sup\{|f(x)|: x\in M\} \quad (f\in A^{*})
\end{align*}
The set of all BSE- functions is expressed by $C_{BSE}(\Delta( A))$. 
The BSE- seminorm of $\sigma  \in C_{BSE}(\Delta(A))$, $q_{l}(\sigma)$, is expressed as:
$$
q_{l}(\sigma)=sup \{|\sum_{i=1}^{n}c_{i}\sigma(\varphi_{i})|:P_{M_{l}}(\sum_{i=1}^{n}c_{i}\varphi_{i})\leq 1 ,~ \varphi_{i}\in\Delta(A), c_{i}\in\mathbb C, n\in\mathbb N\}
$$
where
$$
M_{l}:=\{a\in A: p_{l}(a)\leq 1\}
$$
It was shown that  $(C_{BSE}(\Delta(A)),q_{l})$ is a 
semisimple commutative Frechet subalgebra of $C_{b}(\Delta(A))$.
It is easy to prove that 
$$
q_{l}(\sigma)= inf\{\beta_{M}| \mid\sum_{j=1}^{n}c_{j}\sigma(\varphi_{j})\mid\leq \beta_{M}P_{M}(\sum_{j=1}^{n}c_{j}\varphi_{j}),c_{j}\in\mathbb C, \varphi_{j}\in\Delta(A)\}
$$
It is obvious that if $x\in A$ then $\hat{x}\in C_{BSE}(\Delta(A))$ and 
$q_{l}(\hat{x})\leq p_{l}(x)$, where $\hat{x}(\varphi)= \varphi(x)$ for all $\varphi\in\Delta(A)$. The set $M(A)$ with the strong operator topology, is an unital commutative locally convex algebra.
It was shown that  for each $T\in M(A)$ there exists a unique bounded  continuous function $\widehat{T}$ on $\Delta(A)$ expressed as: 
$$\varphi(Tx)=\widehat{T}(\varphi)\varphi(x),$$
for all $x\in A$ and $\varphi\in \Delta(A)$.
By setting $\{\widehat{T}: T\in M(A)\}$, the $\widehat{M(A)}$ is yield.
A commutative  Frechet algebra $A$ is called BSE- Frechet- algebra  if it meets the following condition:
$$\widehat{M(A)}= C_{BSE}(\Delta(A)).$$
 A bounded net ${\{e_{\beta}\}}$ in $A$ is named a bounded $\Delta $- weak  approximate identity for $A$ if $\varphi(ae_{\beta})\rightarrow \varphi(a)$ for all $\varphi\in\Delta(A)$ and $a\in A$, equivalently $\varphi(e_{\beta})\rightarrow 1$.
\begin{pro}
Let $(A, p_{l})$ be a commutative semisimple Frechet algebra. Then
 $\sigma\in C_{BSE}(\Delta(A))$ if and only if there exists a bounded net $\{x_{\lambda}\}$ in $A$ with 
$$
lim\hat{x_{\lambda}}(\varphi)= \sigma(\varphi)
$$
for all $\varphi\in\Delta(A)$.
\end{pro}\label{prs}
\begin{thm}\label{t.app}
Let $(A, p_{l})$ be a commutative semisimple Frechet algebra. Then
 $A$ has a bounded $\Delta $- weak  approximate identity if and only if
$$\widehat{M(A)}\subseteq C_{BSE}(\Delta(A))$$
\end{thm}


\section{Some basic properties of  vector- valued Frechet- Lipschitz   algebra }\label{2}\setcounter{equation}{0}
       
           The basic terminologies and the related information on vector-valued  Frechet- Lipchitz algebras are reviewed.
 In the sequel, some primary, basic results, and properties related to them are proved. 

Throughout this section, $(X,d)$ is a metric space with at least two elements and $(A,p_{l})$  is a commutative semisimple   Frechet algebra over the scaler field $\mathbb C$. Let $f: X\rt A$ be a function. Set
\begin{align*}
q_{l,A}(f) = \underset{x\in X}{sup} ~ p_{l}(f(x))
\end{align*}
and
\begin{align*}
  p_{l,A}(f) = \underset{x\ne y}{sup} ~ \frac{p_{l}(f(x)-f(y))}{{d(x,y)}}
\end{align*}
 The set of all functions such  $f: X\rt A$ 
satisfies in the following conditions:\\
i) $q_{l,A}(f)<\infty$, for each $ l\in\mathbb N;$\\
ii)$p_{l,A}(f)<\infty$, for each $ l\in\mathbb N.$\\
is  named the vector-valued Frechet Lipschitz algebra and is expressed by ${\Lip}_{d}(X,A)$.

Put\\
$$r_{l,A}(f)=q_{l,A}(f)+p_{l,A}(f)\qquad (f\in {Lip}_{d}(X,A)) $$
 $C_{b}(X,A)$ is the set of all bounded continuous functions from $X$ into $A$. Let $f\in C_{b}(X,A)$.
 If $f,g \in C_{b}(X,A)$ and $\lambda\in\mathbb C$, define\\
$$(f+g)(x)= f(x)+ g(x)\quad and \quad (\lambda f)(x)= \lambda f(x)\quad (x\in X)$$
It is obvious that $(C_{b}(X,A),q_{l,A})$ is a Frechet space over $\mathbb C$. That 
   $({Lip}_{d}(X,A), r_{l,A})$ is a Frechet subalgebra of  $C_{b}(X,A)$ is in Lemma 3.1 proved in  \cite{Ra}.\\
In this article, $f_{a} :X\rt{A}\quad (a\in A)$ is the constant function on $X$, where $f_{a}(x)= a\quad(x\in X)$.
It is obvious that these functions belong to
 $\Lip(X,A)$ and
 $$r_{l,A}(f_{a})= q_{l, A}(f_{a})= p_{l}(a),$$
for each $a\in A$ and $f\in \Lip_{d}(X,A)$.

Let $(A,p_{l})$ and $(B, q_{l})$ be Frechet algebras. The function $\Phi: (A,p_{l})\to (B,q_{l})$ is called isometric, if 
$$q_{l}(\Phi (a))= p_{l}(a) \quad (a\in A,  l\in \mathbb N)$$
Let  $K_{A}:  A\rt {\Lip}_{d}(X,A)$ such that $a\mapsto f_{a}$. Then $K_{A}$ is a continuous, linear and injective function. Furthermore, $\widehat{K_{A}}: \Delta({\Lip}_{\alpha}(d)) \rt \Delta(A)\cup\{0\}$ is homomorphism and $A$ can be considered as a  closed subalgebra of ${\Lip}_{d}(X,A)$.\\

\begin{pro}\label{lem1}
Let $(X,d)$ be a metric space, $(A,p_{l})$ be a commutative Frechet algebra over the scaler field $\mathbb C$. Then ${\Lip}_{d}(X,A)$ is a semisimple Frechet algebra if and only if $A$  is so.
\end{pro}
\begin{proof}
First, Assume that $A$ is semisimple and  take $f,g\in \Lip_{d}(X,A)$ such that $f\ne g$. So there exists $x_{0}\in X$ such that $f(x_{0})\ne g(x_{0})$. Because $A$ is a semisimple algebra,  there exists $\varphi\in\Delta(A)$ where
$$\varphi(f(x_{0}))\ne\varphi(g(x_{0})).$$
  $x_{0}\otimes\varphi$ defined  by $x_{0}\otimes\varphi(f)=\varphi(f(x_{0}))$;  for $f\in \Lip_{d}(X,A)$.  It is obvious that $x_{0}\otimes\varphi\in \Delta(\Lip_{d}(X,A))$ and
$$x_{0}\otimes\varphi(f)\ne x_{0}\otimes\varphi(g)$$
This implies that $\Lip_{d}(X,A)$ is semisimple.\\
Now assume that $\Lip_{d}(X,A)$ is a semisimple Frechet algebra. Let $a,b\in A$ where $a\ne b$, so $f_{a}\ne f_{b}$. 
Because $\Lip_{d}(X,A)$ is a semisimple,  there exists $\psi\in\Delta( \Lip_{d}(X,A))$ where
$\psi(f_{a})\ne \psi(f_{b})$. Which yield:
$$
\widehat{K_{A}}(\psi)(a)=
\psi(K_{A}(a))=\psi(f_{a})\neq \psi(f_{b}) =\psi(K_{A}(b))= \widehat{K_{A}}(\psi)(b) .
$$
Consequently $\widehat{K_{A}}(\psi)(a)\neq \widehat{K_{A}}(\psi)(b)$ and $\widehat{K_{A}}(\psi)\in \Delta(A)$.
Then $\Delta(A)$ separates the points of $A$, this implies that $A$ is a semisimple Frechet algebra.
\end{proof}
\begin{pro}
Let $(X,d)$ be a metric space, $(A,p_{l})$ be a commutative Frechet algebra over the scaler field $\mathbb C$. Then ${\Lip}_{d}(X,A)$ is without order if and only if $A$  is so.
\end{pro}
\begin{proof}
Let $A$ be without order Frechet algebra and Assume that $f\in \Lip_{d}(X, A)$  be non- zero. So there exists $x_{0}\in X$ where $f(x_{0})\ne 0$. Because $A$ is without order, there exists $b\in A$ where
$$f(x_{0})b\ne0.$$
Which yield:
\begin{align*}
(ff_{b})(x_{0})= f(x_{0})f_{b}(x_{0})= f(x_{0})b\ne 0.
\end{align*}
So $ff_{b}\ne 0$., therefore $\Lip_{d}(X,A)$ is without order.\\
Conversely, assume that ${\Lip}_{d}(X,A)$ be without order and  take $a\in A$ where $a\ne 0$, so $f_{a}\ne 0$. Because ${\Lip}_{d}(X,A)$ is without order, there exists $g\in {\Lip}_{d}(X,A)$ where $f_{a}g\ne 0$.
This follows that there exists $x_{0}\in X$ where $(f_{a}g)(x_{0})\ne 0$, thus $ag(x_{0})\ne 0$ and consequently $A$  is without order.
\end{proof}
\begin{lem}\label{l.ap}
Let $(A, p_{l})$ be a commutative semisimple Frechet algebra and $(X,d)$ be a metric space. If  $\Lip_{d}(X,A)$   has a bounded $\Delta-$weak approximate identity, then $A$   has a bounded $\Delta-$weak approximate identity.
\end{lem}
\begin{proof}
Assume that $\Lip_{d}(X, A)$ has a bounded $\Delta-$ weak approximate identity and $(f_{\beta})$ is a bounded $\Delta-$ weak approximate identity for $\Lip_{d}(X, A)$.
By allowing $\varphi\in\Delta(A)$ the following is yield:
$$\underset{\beta}{lim} ~ \varphi(f_{\beta}(x))= \underset{\beta}{lim} ~ (x\otimes\varphi)(f_{\beta})=1.$$
because $x\otimes\varphi\in\Delta(\Lip_{d}(X, A))$, for each $x\in X$ and $\varphi\in \Delta(A)$, thus, the net $(f_{\beta}(x))$ is a bounded $\Delta-$ weak approximate identity for $A$. This completes the proof.
\end{proof}
\begin{lem}
Let $(A, p_{l})$ and $(B, q_{l})$ be a commutative  Frechet algebra and $(X,d)$ be a metric space.
If $A\cong B$, as two Frechet algebras, then $\Lip_{d}(X,A)\cong \Lip_{d}(X,B)$, These two as Frechet algebras are isometric.
\end{lem}
\begin{proof}
Assume that $\Theta: ~ A\longrightarrow B$ is an isomorphism map. Define 
\begin{align*}
\tilde{\Theta}: ~ \Lip_{d}(X,A)\longrightarrow \Lip_{d}(X,B)
\end{align*}
Where $\tilde{\Theta}(f)(x)= \Theta(f(x))$, for all $x\in X$ and $f\in \Lip_{d}(X,A)$. If $f_{1},f_{2}\in \Lip_{d}(X,A)$ and $f_{1}= f_{2}$, so $f_{1}(x)= f_{2}(x)$ 
for each $x\in X$. Thus $\Theta(f_{1}(x))= \Theta(f_{2}(x))$, then $\tilde{\Theta}(f_{1})=\tilde{\Theta}(f_{2})$. This implies that $\tilde{\Theta}$ is well- defined.
At this stage, $\tilde{\Theta}(f)\in \Lip_{d}(X,B)$, for each $f\in \Lip_{d}(X,A)$ is assesed.
For all $x,y\in X$ with $x\neq y$ the following is yield:
\begin{align*}
\frac{q_{l}(\tilde{\Theta}(f)(x)- \tilde{\Theta}(f)(y))}{{d(x,y)}}
                                                                                                            & =\frac{q_{l}({\Theta}(f(x))- {\Theta}(f(y)))}{{d(x,y)}}\\
                                                                                                            &\leq K \frac{p_{l}(f(x)- f(y))}{{d(x,y)}}    
\end{align*}
This follows that
$$
p_{l, B}(\tilde{\Theta}(f))\leq K p_{l, A}(f)
$$ 
Moreover, for all $x\in X$ Which yield: 
$$
q_{l}(\tilde{\Theta}(f)(x))= q_{l}(\Theta(f(x)))\leq Kp_{l}(f(x))
$$
This implies that $q_{l, B}(\tilde{\Theta}(f))\leq Kq_{l,A}(f)$, which $K$ is an upper bound for $\Theta$. Therefore $\tilde{\Theta}(f)\in \Lip_{d}(X,B)$.

In the sequel, it will be concluded that  $\tilde{\Theta}$ is injective. To that end, 
take $f,g\in \Lip_{d}(X,A)$, such that $\tilde{\Theta}(f)= \tilde{\Theta}(g)$. So $\Theta(f(x))= \Theta(g(x))$,
for all $x\in X$, thus $f(x)= g(x)$ for all $x\in X$, because $\Theta$ is injective. Then $\tilde{\Theta}$ is injective.
It remains to prove that  $\tilde{\Theta}$ is surjective. Assume that $g\in \Lip_{d}(X,B)$ and 
define $f(x)= {\Theta}^{-1}(g(x))$, for all $x\in X$. Which yield:
\begin{align*}
\frac{p_{l}(f(x)-f(y))}{{d(x,y)}} &= \frac{p_{l}({\Theta}^{-1}(g(x))-{\Theta}^{-1}(g(y)))}{{d(x,y)}}\\
                                                             &\leq M \frac{q_{l}({g(x)-g(y))}}{{d(x,y)}}
\end{align*}
This follows that 
$$
p_{l, A}(f)\leq M p_{l, B}(g)
$$
In the same way, It will be concluded that $q_{l, A}(f)\leq M q_{l, B}(g)$, for some $M>0$. At the result $f\in \Lip_{d}(X,A)$ and $\tilde{\Theta}(f)(x)= \Theta(f(x))=\Theta({\Theta}^{-1}(g(x)))= g(x)$ and thus $\tilde{\Theta}(f)= g$.
This completes the proof. 
\end{proof}
\begin{lem}
Let $(X,d)$ be a metric space, $(A,p_{l})$ be a commutative Frechet algebra over the scaler field $\mathbb C$. Assume that $M$ is a  bounded set  in $A$, $x\in X$, $\varphi\in\Delta(A)$,  $ c_{1},\cdots,c_{n}\in\mathbb C$  and  the same number $\varphi_{1},\cdots,\varphi_{n}\in\Delta(A)$, 
then the following is yield:
\begin{align*}
P_{M^{'}}(\sum_{i=1}^{n}c_{i}(x\otimes{\varphi}_{i}))= P_{M}(\sum_{i=1}^{n}c_{i}\varphi_{i}).
\end{align*}
where $M^{'}= \{K_{A}(a)| a\in M\}$.
\end{lem}
\begin{proof}
By allowing
 $ c_{1},\cdots,c_{n}\in\mathbb C$  and  the same number $\varphi_{1},\cdots,\varphi_{n}\in\Delta(A)$, the following is yield:
\begin{align*}
P_{M^{'}}(\sum_{i=1}^{n}c_{i}(x\otimes{\varphi}_{i})) 
&= \sup\left\{\left| \sum_{i=1}^{n} c_{i}(x\otimes\varphi_{i})(f)\right|: f\in M^{'}\right\}\\
                         &= \sup\left\{\left|\sum_{i=1}^{n} c_{i}\varphi_{i}(f(x))\right| :  f(x)\in M\right\}\\
                         &\leq \sup\left\{\left|\sum_{i=1}^{n} c_{i}\varphi_{i}(a)\right| :  a\in M\right\}\\ 
             &= P_{M}(\sum_{i=1}^{n}c_{i}\varphi_{i}).
\end{align*}
For the reverse inclusion,  which yield:
\begin{align*}
 P_{M}(\sum_{i=1}^{n}c_{i}\varphi_{i}) &=  \sup\left\{\left|\sum_{i=1}^{n} c_{i}\varphi_{i}(a)\right| :  a\in M\right\}\\
                                                                                        &=  \sup\left\{\left|\sum_{i=1}^{n} c_{i}\varphi_{i}(f_{a}(x))\right| :  a\in M\right\}\\
                                                                                        &=  \sup\left\{\left|\sum_{i=1}^{n} c_{i}(x\otimes{\varphi}_{i})(f_{a})\right| :  f_{a}\in M^{'}\right\}\\
                                                                                       &\leq \sup\left\{\left| \sum_{i=1}^{n} c_{i}(x\otimes\varphi_{i})(f)\right|: f\in M^{'}\right\}\\
                                                                                       &= P_{M^{'}}(\sum_{i=1}^{n}c_{i}(x\otimes{\varphi}_{i}))
\end{align*}
Consequently,
$$
P_{M^{'}}(\sum_{i=1}^{n}c_{i}(x\otimes{\varphi}_{i}))= P_{M}(\sum_{i=1}^{n}c_{i}\varphi_{i}).
$$
\end{proof}

\setcounter{section}{2}
\section{ Main results  }\label{2}\setcounter{equation}{0} 
              The structure of the BSE functions on $\Delta(\Lip_{d}(X,A))$ is characterization and the
correlations between the BSE property of A and $\Lip_{d}(X,A)$ are assessed.

Let $f\in \Lip_{d}(X,A)$, define $r^{'}_{l, A}(f)= max\{p_{l, A}(f), q_{l, A}(f)\}$. It is obvious that $r^{'}_{d, A}$ is a seminorm on $\Lip_{d}(X,A)$.  
Clearly 
$$(\Lip_{d}(X,A),r_{l, A})\cong ( \Lip_{d}(X,A),r_{l, A}^{'})$$
\begin{pro}\label{gh1}
Let $(X,d)$ be a metric space and  $(A,p_{l})$ be a  commutative semisimple Frechet algebra. Assume that $\Lip_{d}(X,A)$ is a BSE- Frechet- algebra. Then $A$ is so.
\end{pro}
\begin{proof}
Because $\Lip_{d}(X,A)$ is a BSE- algebra, by  referring to Theorem\ref{t.app}, $\Lip_{d}(X,A)$ has a bounded $\Delta-$ weak approximate identity.   Lemma\ref{l.ap} and Theorem \ref{t.app}, implies that
$$\widehat{M(A)}\subseteq C_{BSE}(\Delta(A)).$$
For the reverse inclusion, take $\sigma\in C_{BSE}(\Delta(A))$. There exist a bounded set $M$ in $A$ and a positive real number $\beta_{M}$ where by allowing  $\psi_{1},\cdots,\psi_{n}$ of $\Delta(\Lip_{d}(X,A))$ 
and the same number of complex numbers $c_{1},\cdots,c_{n}$, the following is yield:
\begin{align*}
\mid\sum_{i=1}^{n}c_{i}\sigma o\widehat{K_{A}}(\psi_{i})\mid &= \mid\sum_{i=1}^{n}c_{i}\sigma(\psi_{i}o K_{A})\mid\\
                                                                                             &\leq\beta_{M}P_{M}(\sum_{i=1}^{n}c_{i}(\psi_{i}o K_{A}))\\
                                                                                              &\leq\beta_{M}KP_{M^{'}}(\sum_{i=1}^{n}c_{i}\psi_{i})
\end{align*}
for some $K>0$, where $M^{'}= \{K_{A}(a)| a\in M\}$. It follows that $\sigma o\widehat{K_{A}}\in C_{BSE}(\Delta(\Lip_{d}(X,A)))$.
By applying the BSE- property of  $\Lip_{d}(X,A)$,  there exists $T\in M(\Lip_{d}(X,A))$ where $\hat{T}=\sigma o\widehat{K_{A}}$. Now define $T^{'}\in M(A)$ as follows:
$$T^{'}(a)=T(K_{A}(a))(x_{0}),\quad (a\in A)$$
where $x_{0}\in X$ is an arbitrary member of $X$. If $a_{1}, a_{2}\in A$;
\begin{align*}
 T^{'}(a_{1}a_{2})= T(K_{A}(a_{1}a_{2}))(x_{0}) &= T(K_{A}(a_{1}) K_{A}(a_{2}))(x_{0})\\
                                                                                &= (T(K_{A}(a_{1})) K_{A}(a_{2}))(x_{0})\\
                                                                                 &= T(K_{A}(a_{1}))(x_{0}) K_{A}(a_{2})(x_{0})= T^{'}(a_{1}).a_{2}
\end{align*}
Hence $T^{'}\in M(A)$.
Let $\varphi\in \Delta(A)$;  It is easy to prove that $\hat{K_{A}}(x_{0}\otimes\varphi)=\varphi$ and the following is yield:
\begin{align*}
\hat{ T^{'}}(\varphi) =  \frac{\varphi(T^{'}(a))}{\varphi(a)}&= \frac{\varphi(T(K_{A}(a))(x_{0}))}{\varphi(a)}\\
                                                                                                  &=  \frac{(x_{0}\otimes\varphi)(T(K_{A}(a)))}{\varphi(f_{a}(x_{0}))}\\
                                                                                                  &=  \frac{(x_{0}\otimes\varphi)(T(K_{A}(a)))}{(x_{0}\otimes\varphi)(K_{A}(a))}\\
                                                                                                  &= \hat{T}(x_{0}\otimes\varphi)\\
                                                                                                  &= \sigma o\widehat{K_{A}}(x_{0}\otimes\varphi)\\
                                                                                                  &=\sigma(\widehat{K_{A}}(x_{0}\otimes\varphi))\\
                                                                                                   &= \sigma(\varphi)
\end{align*}
Therefore $\hat{T^{'}}= \sigma$ and consequently $C_{BSE}(\Delta(A))\subseteq\widehat{M((A)}$. Thus  $A$ is a Frechet- BSE- algebra.
\end{proof}
The correlation between the $C_{\textup{BSE}}(\Delta(\Lip_{d}(X,A)))$ and
 $\Lip_{d}(X,C_{\textup{BSE}}(\Delta(A)))$ is assessed as follows:
\begin{thm}\label{thc}
Let $(X,d)$ be a metric space,  $(A,p_{l})$ be a  commutative semisimple Frechet algebra.  Then
$C_{BSE}(\Delta(\Lip_{d}(X,A)))$ can be embedded in $\Lip_{d}(X,C_{BSE}(\Delta(A))$, These two as Frechet algebras are isometric;
\end{thm}
\begin{proof}
Let
\begin{align*}
\phi: ~ C_{BSE}(\Delta(\Lip_{d}(X,A)))\to \Lip_{d}(X,C_{BSE}(\Delta(A))) ,    
\end{align*}
defined by 
\begin{align*}
\phi(\Sigma)= \phi_{\Sigma}    \quad (\Sigma\in C_{BSE}(\Delta(\Lip_{d}(X,A))),
\end{align*}
Where
\begin{align*}
  \phi_{\Sigma}(x)(\varphi)= \Sigma(x\otimes\varphi),\quad (x\in X, ~\varphi\in\Delta(A))
\end{align*}
Assume that $\Sigma_{1},\Sigma_{2}\in C_{BSE}(\Delta(\Lip_{d}(X,A)))$, where $\Sigma_{1} =\Sigma_{2}$.\\
So $\Sigma_{1}(x\otimes\varphi)= \Sigma_{2}(x\otimes\varphi)$, at the result $\phi_{\Sigma_{1}}(x)(\varphi)= \phi_{\Sigma_{2}}(x)(\varphi)$, for all $x\in X$ and $\varphi\in\Delta(A)$. Then $\phi_{\Sigma_{1}}= \phi_{\Sigma_{2}}$ and, therefore, $\phi$ is well defined.
It is obvious that $\phi$ is linear. Let $\Sigma_{1},\Sigma_{2}\in C_{BSE}(\Delta(\Lip_{d}(X,A)))$, so 
$\phi(\Sigma_{1}.\Sigma_{2})= \phi_{\Sigma_{1}.\Sigma_{2}}$, By allowing $x\in X$ and $\varphi\in\Delta(A)$,
the following is yield:
\begin{align*}
\phi_{\Sigma_{1}.\Sigma_{2}}(x)(\varphi) &= \Sigma_{1}.\Sigma_{2}(x\otimes \varphi)\\
                                                                &= \Sigma_{1}(x\otimes \varphi). \Sigma_{2}(x\otimes \varphi)\\
                                                                &= \phi_{\Sigma_{1}}(x)(\varphi). \phi_{\Sigma_{2}}(x)(\varphi)
\end{align*}
 thus $\phi_{\Sigma_{1}.\Sigma_{2}}(x)= \phi_{\Sigma_{1}}(x). \phi_{\Sigma_{2}}(x)$, 
so $\phi_{\Sigma_{1}.\Sigma_{2}}= \phi_{\Sigma_{1}}. \phi_{\Sigma_{2}}$. Then $\phi$ is homomorphism.
First of all,  
$\phi_{\Sigma}(x)\in C_{BSE}(\Delta(A))$, for each $x\in X$ and $\Sigma\in C_{BSE}(\Delta(\Lip_{d}(X,A))$ is assessed. In fact, Since $\Sigma\in C_{BSE}(\Delta(\Lip_{d}(X,A))$, so there exists a bounded set $M$ in $\Lip_{d}(X,A)$  such that for every complex number
 $ c_{1},\cdots,c_{n}$  and  the same number $\varphi_{1},\cdots,\varphi_{n}\in\Delta(A)$, we have
\begin{align*}
{P}_{{M}}(\sum_{i=1}^{n}c_{i}(x\otimes{\varphi}_{i})) &= sup\{~ | \sum_{i=1}^{n} c_{i}(x\otimes\varphi_{i})(f)|: ~ f\in M\}\\
                                                                                                            &= sup\{ ~ |\sum_{i=1}^{n} c_{i}\varphi_{i}(f(x))| : ~ f\in M\}\\
                                                                                                            &\leq sup \{ ~| \sum_{i=1}^{n} c_{i}\varphi_{i}(a)| : ~ a\in M^{'}\}\\ 
                                                                                                             &= P_{M^{'}}(\sum_{i=1}^{n}c_{i}\varphi_{i}).
\end{align*}
Where $M^{'}:= \hat{x}(M)$.
This implies that
\begin{align*}
|\sum_{i=1}^{n}c_{i}\phi_{\Sigma}(x)(\varphi_{i})| &= |\sum_{i=1}^{n}c_{i}\Sigma(x\otimes\varphi_{i})|\\
                                                                                     &\leq q_{l}(\Sigma){P}_{{M}}(\sum_{i=1}^{n}c_{i}x\otimes{\varphi}_{i})\\
                                                                                     &= q_{l}(\Sigma)P_{M^{'}}(\sum_{i=1}^{n}c_{i}\varphi_{i})
\end{align*}
Hence $\phi_{\Sigma}(x)\in C_{BSE}(\Delta(A))$ and
 $q_{l}(\Sigma)\leq q_{l}(\phi_{\Sigma}(x))$, for each $x\in X$, since  $q_{l}(\Sigma)= inf\{\beta_{M} | \mid\sum_{j=1}^{n}c_{j}\Sigma(\varphi_{j})\mid\leq \beta_{M}P_{M}(\sum_{j=1}^{n}c_{j}\varphi_{j})\}$  In the other hand
\begin{align*}
q_{l}(\phi_{\Sigma}(x)) &=sup \left\{ ~\left|\sum_{i=1}^{n}c_{i}(\phi_{\Sigma}(x))(\varphi_{i})\right|: ~ P_{M_{l}}(\sum_{i=1}^{n}c_{i}\varphi_{i})\leq 1,~\varphi_{i}\in\Delta(A)\right\}\\
                                      &=sup \left\{\left|\sum_{i=1}^{n}c_{i}\Sigma(x\otimes\varphi_{i})\right|: ~ P_{M_{l}}(\sum_{i=1}^{n}c_{i}(x\otimes\varphi_{i}))\leq 1 ,~ \varphi_{i}\in\Delta(A)\right\}\\
                                       &\leq q_{l}(\Sigma)
\end{align*}
Therefore 
$$q_{l}(\phi_{\Sigma}(x))\leq q_{l}(\Sigma).$$
 Consequently, for all $x\in X, \Sigma\in C_{BSE}(\Delta(\Lip_{d}(X, A)))$ and $l\in \mathbb N$, we have
\begin{align}\label{1}
 q_{l}(\phi_{\Sigma}(x))= q_{l}(\Sigma).
\end{align}
Note that
\begin{align*}
q_{l}(\phi_{\Sigma}(x)-\phi_{\Sigma}(y)) &=~ sup\{|\sum_{i=1}^{n}c_{i}(\phi_{\Sigma}(x)-\phi_{\Sigma}(y))(\varphi_{i})|: ~ P_{M_{l}}(\sum_{i=1}^{n}c_{i}\varphi_{i})\leq 1, ~ \varphi_{i}\in\Delta(A)\}\\
                                                                 &= ~ sup\{|\sum_{i=1}^{n}c_{i}(\Sigma(x\otimes\varphi_{i})-\Sigma(y\otimes\varphi_{i}))|: ~ P_{M_{l}}(\sum_{i=1}^{n}c_{i}\varphi_{i})\leq 1, ~\varphi_{i}\in\Delta(A)\}\\
                                                                     &\leq q_{l}(\Sigma)sup\{P_{M_{l}}(\sum_{i=1}^{n}c_{i}(x\otimes\varphi_{i}-y\otimes\varphi_{i})) : ~ P_{M_{l}}(\sum_{i=1}^{n}c_{i}\varphi_{i})\leq 1, ~ \varphi_{i}\in\Delta(A)\}.                                                                                                                       
\end{align*}
This follows that
\begin{align*}
p_{l,C_{BSE}(\Delta(A))}(\phi_{\Sigma}) &= ~ sup\{\frac{q_{l}(\phi_{\Sigma}(x)-\phi_{\Sigma}(y))}{d(x,y)}: x\ne y\} \\
                                                &\leq q_{l}(\Sigma)sup ~sup\{\{\frac {P_{M_{l}}(\sum_{i=1}^{n}c_{i}(x\otimes\varphi_{i}-y\otimes\varphi_{i}))}{{d(x,y)}}|:\\
                                                 & P_{M_{l}}(\sum_{i=1}^{n}c_{i}\varphi_{i}) \leq 1, ~ \varphi_{i}\in\Delta(A)\}: x\ne y\}\\      
                                                &\leq q_{l}(\Sigma)         
\end{align*}
Therefore $p_{l,C_{BSE}(\Delta(A))}(\phi_{\Sigma})\leq q_{l}(\Sigma)$. Also
$$
q_{l,C_{BSE}(\Delta(A))}(\phi_{\Sigma})= sup\{q_{l}(\phi_{\Sigma})(x)): x\in X\}
$$
 and by using relation (\ref{1}),  
$q_{l}(\phi_{\Sigma}(x))=q_{l}(\Sigma)$
 and so $q_{l,C_{BSE}(\Delta(A))}(\phi_{\Sigma})= q_{l}(\Sigma)$. This show that
 $$r_{l,C_{BSE}(\Delta(A))}^{'}(\phi_{\Sigma})= max\{p_{l,C_{BSE}(\Delta(A))}(\phi_{\Sigma}),q_{l,C_{BSE}(\Delta(A))}(\phi_{\Sigma})\}=q_{l}(\Sigma)$$
Which implies that $\phi$ is isometry.  This completes the proof.
\end{proof}
Let $T\in M(A)$. Define
$$q_{l}^{'}(T)= sup\{ p_{l}(T(a)): a\in A,~ p_{l}(a)\leq 1\}.$$
It is obvious that $q_{l}^{'}$ is a seminorm on $M(A)$.
In the following theorem, it will be shown that $Lip_{\alpha}(X,A)$ can be embedded in
$M(\Lip_{d}(X,A))$, isometrically as two locally convex algebras which are isometric.
\begin{thm}\label{thM}
Let $(X,d)$ be a metric space and  $(A,p_{l})$ be a  commutative semisimple Frechet algebra.
 Then 
$$\Lip_{d}(X, M(A))\underset{\rightarrow}{\subseteq} M(\Lip_{d}(X,A)),$$
As two locally convex algebras which are isometric.
\end{thm}
\begin{proof}
Let
\begin{align*}
\phi ~:  \Lip_{d}(X,M(A))\to M(\Lip_{d}(X,A)) 
\end{align*}
Where
\begin{align*}
  \phi(F)= \phi_{F}\quad      ( F\in \Lip_{d}(X,M(A))).
\end{align*}
Defined by
\begin{align*}
& \phi_{F}(g)= F\odot g  \quad   & ( g\in \Lip_{d}(X, A))\\
& F\odot g(x)= F(x)(g(x)) \quad & ( x\in X).
\end{align*}
It will be concluded that $\phi$ is an isomorphism map.
Assume that $F_{1}, F_{2}\in \Lip_{d}(X,M(A))$ where $F_{1}= F_{2}$, so $F_{1}(x)= F_{2}(x)$,
for each $x\in X$. Thus  $F_{1}(x)(g(x))= F_{2}(x)(g(x))$, for all $g\in \Lip_{d}(X, A)$, then
$F_{1}\odot g(x)= F_{2}\odot g(x)$, for all $x\in X$ and $g\in \Lip_{d}(X, A)$.
Therefore $\phi(F_{1})=\phi(F_{2})$ and so $\phi$ is well-defined.\\
1) $\phi_{F}$ is a continuous linear multiplier on $Lip_{\alpha}(X, A)$ is assessed in the following:\\
Assume that $g_{1}, g_{2}\in \Lip_{d}(X,M(A))$, so 
$$
\phi_{F}(g_{1}g_{2})= F\odot g_{1}g_{2}.
$$
for any $x\in X$, the following is yield:
\begin{align*}
F\odot g_{1}g_{2}(x) &= F(x)(g_{1}g_{2}(x))\\ 
                                   &= F(x)(g_{1}(x)g_{2}(x)) \\
                                   &= g_{1}(x)F(x)(g_{2}(x))\\ 
                                   &= g_{1}(x)F\odot g_{2}(x)  
\end{align*}
This implies that
$$
 F\odot g_{1}g_{2}= g_{1}.(F\odot g_{2}).
$$
Then
$
\phi_{F}(g_{1}g_{2})= g_{1}.\phi_{F}(g_{2})
$, for all $g_{1}, g_{2}\in \Lip_{d}(X,M(A))$, at the result $\phi_{F}$ is a multiplier.
It is obvious that for each $F\in \Lip_{d}(X,M(A))$,  the map $\phi_{F}$ is linear.
In the sequel,  that $\phi_{F}$ is a continuous map will be proved. Let $(g_{n})$ be a sequence in $\Lip_{d}(X, A)$  converges to $g\in \Lip_{d}(X,A)$. Then $r_{l,A}(g_{n})\to r_{l,A}(g)$, thus $q_{l,A}(g_{n})\to q_{l,A}(g)$ and\\
 so $p_{l,A}(g_{n})\to p_{l,A}(g)$.
Which yield:
\begin{align*}
q_{l,A}(F\odot g_{n}-F\odot g) &= {sup}\left\{p_{l}(F(x)(g_{n}(x)-g(x))) ~: ~ x\in X\right\}\\
                                                          &\leq K ~ {sup}\left\{p_{l}(g_{n}(x)-g(x))~: ~ x\in X\right\}\\
                                                           &= Kq_{l,A}(g_{n}-g)\to 0,
\end{align*}
Which
\begin{align*}
K:= q_{l,M(A)}(F)=sup\left\{q_{l}^{'}(F(x))| x\in X\right\}
\end{align*}
Note that
\begin{align*}
q_{l}^{'}(F(x))= {sup}\left\{p_{l}(F(x))(a)| a\in A, p_{l}(a)\leq 1\right\}.
\end{align*}
This shows that
\begin{align*}
q_{l,A}(F\odot g_{n}-F\odot g)\to 0.
\end{align*}
 Also:
\begin{align*}
p_{l,A}(F\odot g_{n}-F\odot g) &={sup}\left\{\frac{p_{l}(F(x)(g_{n}(x)-g(x))-F(y)(g_{n}(y)-g(y)))}{d(x,y)}~: ~ {x\ne y}\right\}\\
                                                         &\leq {sup}\left\{\frac{p_{l}\left(F(x)(g_{n}(x)-g(x))-(g_{n}(y)-g(y))\right)}{d(x,y)} ~: ~ {x\ne y}\right\}\\
                                                         &+ {sup}\left\{\frac{p_{l}(F(x)-F(y))}{d(x,y)})(g_{n}(y)-g(y): ~ {x\ne y}\right\}\\
                                                         &\leq q_{l,M(A)}(F)p_{l,A}(g_{n}-g)+p_{l,M(A)}(F)q_{l,A}(g_{n}-g)\to 0
\end{align*}
Which yield:
\begin{equation*}
p_{l,A}(F\odot g_{n}-F\odot g)\to 0.
\end{equation*}
Therefore $r_{l,A}(F\odot g_{n}-F\odot g)\to 0$. Hence $\phi_{F}$ is a continuous map.
This follows that $\phi_{F}\in M(\Lip_{d}(X,A))$.\\
2) In following, it will be concluded that $\phi$ is an isomorphism map.\\
It is obvious that $\phi$ is a linear map. Assume that $(F_{n})$ is a sequence  in $\Lip_{d}(X,M(A))$ which converges to some $F\in \Lip_{d}(X,M(A))$, so
$p_{l,M(A)}(F_{n}-F)\to 0$ and $q_{l,M(A)}(F_{n}-F)\to 0$. yielding the following:
\begin{align*} 
q_{l}^{'}(\phi(F_{n})-\phi(F)) &= sup\left\{r_{\alpha,A}(F_{n}\odot g- F\odot g): {g\in \Lip_{d}(X,A)}\right\}\\
\end{align*}
Hence
\begin{align*}
q_{l,A}(F_{n}\odot g- F\odot g) &=sup\left\{p_{l}((F_{n}(x)-F(x))g(x)): {x\in X}\right\} \\
                                                    &\leq {sup}{sup}\left\{p_{l}((F_{n}-F)(x)(a)): a\in A, {x\in X}\right\}\\
                                                    &={sup}\left\{q_{l}^{'}((F_{n}-F)(x)): {x\in X}\right\}= q_{l,M(A)}(F_{n}-F)\to 0
\end{align*}
Thus
\begin{align*}
 q_{l,A}(F_{n}\odot g- F\odot g)\to 0.
\end{align*}
Moreover
\begin{align*} 
p_{l,A}(F_{n}\odot g- F\odot g) &={sup}\left\{\frac{p_{l}((F_{n}-F)(x)(g(x)))- ((F_{n}-F)(y)(g(y)))}{d(x,y)}: {x\ne y}\right\}\\
                                                            &\leq {sup}\left\{\frac{p_{l}(((F_{n}-F)(x)-(F_{n}-F)(y))(g(x))}{d(x,y)} : {x\ne y}\right\}\\
                                                            &+{sup}\left\{\frac{(p_{l}(F_{n}-F)(y)(g(x) -g(y))}{d(x,y)} : {x\ne y}\right\}\\ 
                                                            &\leq p_{l,M(A)}(F_{n}-F).q_{l,A}(g)+q_{l,M(A)}(F_{n}-F).p_{l,A}(g)\rightarrow  0,
\end{align*}
and so
\begin{align*}
 p_{l,A}(F_{n}\odot g- F\odot g)\to 0.
\end{align*}
 It follows that $r_{l,A}(\phi_{F_{n}}(g)-\phi_{F}(g))\to  0$, for all $g\in \Lip_{d}(X,A)$. Consequently\\
 $q_{l}^{''}(\phi(F_{n})-\phi(F))\to 0$. Therefore $\phi$
is a continuous map.\\ 
In this stage, that $\phi$ is injective is assessed. Assume that  $F\in \Lip_{d}(X,M(A))$ where
 $\phi_{F}= \phi(F)= 0$. If $F\ne0$, then there exists $x_{0}\in X$ where $F(x_{0})\ne 0$, so there exists $a_{0}\in A$ where $(F(x_{0}))(a_{0})\ne 0$. Put $g=f_{a_{0}}$, thus
$g\in \Lip_{d}(X,A)$ and
$$F\odot g(x_{0})= F(x_{0})(g(x_{0}))=(F(x_{0}))(a_{0}))\ne 0$$
i.e.
$$
\phi_{F}(g)(x_{0})= F\odot g(x_{0})\ne 0
$$
This is a contradiction. Therefore $\phi$ is injective.

\end{proof}
At this stage, based on the established prerequisite the primary Theorem,
is expressed as follows:
\begin{thm}\label{thna}
Let $(X,d)$ be a metric space and  $(A,p_{l})$ be a  commutative semisimple Frechet algebra.
 Then  $\Lip_{d}(X, A)$  is a Frechet- BSE-algebra if and only if $A$ is a Frechet- BSE algebra. 
 Then\\ 
 1)
If 
 $\Lip_{d}(X, A)$  is a   Frechet- BSE-algebra, then $A$ is a  Frechet- BSE-algebra. \\
2) If $A$ is unital and  Frechet- BSE-algebra, then $\Lip_{d}(X, A)$ is a  Frechet- BSE-algebra.

\end{thm}
\begin{proof}
1) If   $\Lip_{d}(X,A)$ is a Frechet- BSE-algebra, then by Proposition\ref{gh1}, $A$ is a Frechet- BSE-algebra.\\
2) 
Assume that $A$ is a   BSE-algebra.
Since $A$ is semisimple, then by applying Proposition\ref{lem1}, $\Lip_{d}(X, A)$ is semisimple. By applying Proposition\ref{gh1}  imply that 
\begin{align*}
{(M(\Lip_{d}(X,A)\widehat{)}}\subseteq C_{\textup{BSE}}(\Delta(\Lip_{d}(X,A)).
\end{align*}

For the reverse inclusion, 
 according to Theorem \ref{thc} and Theorem \ref{thM}, the following is yield:
\begin{align*}
C_{\textup{BSE}}(\Delta(\Lip_{d}(X,A)) &\underset{\rightarrow}{\subseteq} \Lip_{d}(X,C_{\textup{BSE}}(\Delta(A)))\\
                                                        &\cong \Lip_{d}(X,\widehat{M(A)})\\
                                                       &= \Lip_{d}(X,M(A))\\
                                                      &\underset{\rightarrow}{\subseteq} M(\Lip_{d}(X,A)\\
                                                     &\cong {(M(\Lip_{d}(X,A)\widehat{)}}.
\end{align*}
Thus
\begin{align*}
C_{\textup{BSE}}(\Delta(\Lip_{d}(X,A))\cong (M(\Lip_{d}(X,A))\widehat{)}.
\end{align*}
\end{proof}
Because every commutative  $C^{*}$- Banach algebra is BSE algebra, \cite{E7}, so by using Theorem\ref{thna}, 
the following example is immediate: 
\begin{ex}
Let $(X,d)$ be a metric space and $A$ be a commutative  $C^{*}$- Banach algebra. Then $\Lip_{d}(X,A)$ is a BSE- Frechet algebra
\end{ex}


\end{document}